\documentclass[12pt]{article}
\usepackage[leqno]{amsmath}
\usepackage{amsthm,amscd}
\usepackage{amsfonts} 
\usepackage{amssymb} 
\newtheorem{theorem}{Theorem}[section]

\newtheorem{example}[theorem]{Example}
\newtheorem{lemma}[theorem]{Lemma}

\newtheorem{corollary}[theorem]{Corollary}

\newcommand{\A}{{\mathcal A}}
\newcommand{\B}{{\mathcal B}}

\newcommand{\bfx}{{\mathbf x}}
\newcommand{\bfy}{{\mathbf y}}

\newcommand{\bbP}{{\mathbb P}}

\newcommand{\cR}{{\mathcal R}}
\newcommand{\bbK}{{\mathbb K}}
\newcommand{\bbF}{{\mathbb F}}
\newcommand{\bbR}{{\mathbb R}}

\newcommand{\bbS}{{\mathbb S}}

\newcommand{\be}{{\mathbf e}}
\newcommand{\bfz}{{\mathbf z}}

\newcommand{\im}{\operatorname{im}}

\newcommand{\ch}{{\rm \bf Ch}}

\newcommand{\ra}{\rightarrow}
\newcommand{\ints}{\mathbb Z}

\usepackage{graphicx}

\begin{document}

\title{Ranking Patterns of the Unfolding Model and Arrangements}

\author{Hidehiko Kamiya
\\
{\it Okayama University, Okayama, Japan} \\
Peter Orlik \\
{\it University of Wisconsin, Madison, WI, USA} \\
Akimichi Takemura \\
{\it University of Tokyo, Tokyo, Japan} \\
and \\
Hiroaki Terao \\
{\it Tokyo Metropolitan University, Tokyo, Japan} \\
}

\date{April, 2004}

\maketitle

\begin{abstract}

In the unidimensional unfolding model, given $m$ objects in
general position there arise $1+m(m-1)/2$ rankings. The set of
rankings is called the ranking pattern of the $m$ given objects.
By changing these $m$ objects, we can generate various ranking
patterns. It is natural to ask how many ranking patterns can be
generated and what is the probability of each ranking pattern when
the objects are randomly chosen? These problems are studied by
introducing a new type of arrangement called mid-hyperplane
arrangement and by counting cells in its complement.



\end{abstract}

\maketitle

\bigskip




\section{Introduction}

Various models have been developed for the analysis of ranking
data. These include Thurstonian models, distance-based models,
paired and higher-order comparison models, ANOVA-type loglinear
models, multistage models and unfolding models, to name only a
few. These models give a description of the ranking process and/or
the population of rankers. For a comprehensive treatment of the
methods for analyzing and modelling ranking data, see the
excellent book by Marden~\cite{mar}.

The unfolding model was devised by Coombs~\cite{coo1, coo2, coo3}
for the analysis of ranking data based on preferential choice
behavior. According to De Soete, Feger and Klauer~\cite[p.1]{dfk},
``Historically, two of the most important contributions to
psychological choice modelling are undoubtedly
Thurstone's~\cite{thu} Law of Comparative Judgment and
Coombs'~\cite{coo1, coo3} unfolding theory.'' This model has been
widely used in practice in many fields beyond psychology:
sociology, marketing science, voting theory, etc. In addition, the
same mathematical structure can be found in Voronoi diagrams
(Okabe, Boots, Sugihara and Chiu~\cite{obs}), spatial competition
models in urban economics (Hotelling~\cite{hot}, Eaton and
Lipsey~\cite{eal1, eal2}) and multiple discriminant analysis
(Kamiya and Takemura~\cite{kat}).

According to the unidimensional unfolding model, preferential
choice is made in the following manner: all individuals evaluate
$m$ objects based on the objects' single common attribute. Each
object is represented by a real number expressing the level of
this attribute $x_i, \ i=1,2,\ldots, m,$ or a point on the real
line $\bbR$ (the ``unidimensional underlying continuum''). At the
same time, each individual is also represented by a point $y \in
\bbR$ on the same line. The point $y$ is considered the
individual's favorite and is called his/her ideal point. In this
model, the real line $\bbR$ containing both individuals and
objects is thought of as the psychological space and is called the
{\it joint scale} or the J scale. Here we identify individuals and
objects with their corresponding points. The model assumes that
individual $y$ ranks the $m$ objects $x_1, x_2, \ldots, x_m$
according to their distances from $y,$ i.e., individual $y$
prefers $x_i$ to $x_j$ iff $|y-x_i|<|y-x_j|.$ Rankings generated
by individuals in this way are sometimes called individual scales
or I scales.

We say that the $m$ points representing the objects are in general
position if they  and their midpoints are all distinct. Further,
we do not consider partial rankings or ties in this paper, so we
treat only those individuals whose ideal points do not coincide
with any midpoint of two objects.

Let  $x_1,\ldots,x_m$  be  $m$ objects which satisfy these
assumptions. By varying the location of the ideal point $y$
throughout $\bbR$ except the midpoints, we can account for
$\binom{m}{2}+1$ kinds of rankings of $x_1,x_2,\ldots,x_m$.  The
significance of using this model lies here: there are $m!$
potential rankings, but the psychological structure restricts the
variety of rankings that can actually occur. These
$\binom{m}{2}+1$ rankings are called the {\it admissible rankings}
of $\bfx=(x_1,\ldots,x_m).$

The unidimensional unfolding model was extended to the
multidimensional case by Bennett and Hays~\cite{beh} and Hays and
Bennett~\cite{hab}. As the dimension $n$ of the psychological
space gets large, the number of admissible rankings accounted for
by this model increases, hence more rankings can be explained.
This means that the psychological structure becomes looser as $n$
increases. In fact, when $n \ge m-1,$ all $m!$ rankings are
admissible, and the model in this case is not interesting at all.
Thus finding an appropriate dimension is important in the actual
analysis of ranking data. In this paper, we restrict our attention
to the unidimensional unfolding model.

For a given set of $m$ objects represented by $\bfx=(x_1, \ldots,
x_m)\in \bbR^n,$ we call the set of $\binom{m}{2}+1$ admissible
rankings of $\bfx$ the {\it ranking pattern} of $\bfx.$ By
considering different attributes, we can get different sets of $m$
real numbers $x_1,x_2,\ldots,x_m$ for the same $m$ objects, and
thus obtain different ranking patterns. Examination of the
collected sample of ranking data can tell us what aspect of the
$m$ objects determines the present individuals' preferential
choice behavior towards these $m$ objects, thereby enabling some
inference about the latent structure.

It is generally impossible to explain all  ranking data by
considering any single attribute. Van
Blokland-Vogelesang~\cite{vbl} introduces an error structure into
the unfolding model and makes it a probabilistic model. The error
structure in her model is an extension of the Feigin and
Cohen~\cite{fec} model. Other types of error structures have also
been studied by other authors (Brady~\cite{bra}, B\"{o}ckenholt
and Gaul~\cite{bog}, De Soete, Carroll and DeSarbo~\cite{dcd}).
Moreover, for a set of ranking data which is not completely
compatible with any joint scale, van Blokland-Vogelesang~\cite{vbl} 
proposes a method for finding the ``best'' joint scale based on
Kendall's $\tau$ distance.

The arguments so far imply that it is important to know the
variety of ranking patterns generated by the unfolding model. The
significance of this problem can also be understood in the context
of voting theory or social choice theory. It is well known
(Coombs~\cite{coo3}, Luce and Raiffa~\cite{lur}) that the
unfolding model avoids voting cycles by restricting  the possible
rankings.

Suppose three individuals $A, B$ and $C$ rank three objects
labelled $1, 2$ and $3$ as $(123), (231)$ and $(312),$
respectively, where objects $1, 2$ and $3$ are listed in order
from best to worst in the expression $(i_1 i_2 i_3).$ Here two
individuals $A$ and $C$ prefer $1$ to $2,$ while $B$ prefers $2$
to $1,$ so by  simple majority rule, $1$ is preferred to $2$ as a
collective preference. In the same way, the simple majority rule
yields the collective preference that $2$ is preferred to $3$ and
that $3$ is preferred to $1,$ entailing intransitivity called a
voting cycle.

But if the individuals' preferences are limited to those determined
by the unfolding model, we can see that the collective preference
by  simple majority rule coincides with the median individual's
preference and thus in particular produces no voting cycles. Here
the median individual means the individual $M \in \{ A, B, C \}$
whose ideal point $y_M \in \bbR$ is the median of the individuals'
ideal points $y_A, y_B, y_C \in \bbR.$ The same holds true for any
odd number of individuals and any number $m \ge 3$ of objects.
Thus it is crucial to clarify how much restriction the unfolding
model imposes on individuals' possible preferences.

We show next that it suffices to study the case where the $m$
objects $x_1,\ldots,x_m$ are ordered as $x_1<\cdots<x_m.$ Consider
two sets of $m$ objects $\bfx=(x_1,\ldots,x_m)$ and
$\bfx'=(x'_1,\ldots,x'_m).$ If the rank orders of their midpoints
from left to right on $\bbR$ are the same or the reverse of each
other, then $\bfx$ and $\bfx'$ produce the same ranking pattern.
Conversely, if $\bfx$ and $\bfx'$ induce the same ranking pattern,
their midpoint orders are the same or the reverse of each other.
These facts can be confirmed by using  results of Kamiya and
Takemura~\cite{kat}. Here we agree to say that two rank orders of
midpoints (or objects) are essentially different if one is
different from the other as well as from the reverse of the other.
These arguments imply that there is a one-to-one correspondence
between the set of ranking patterns and the set of essentially
different rank orders of midpoints of the objects. Since the rank
order of objects is completely determined by the rank order of
their midpoints, two sets of $m$ objects having essentially
different rank orders give rise to essentially different rank
orders of their midpoints and thus different ranking patterns.

On the other hand, it is obvious that for any permutation $(i_1
\ldots i_m)$ of $\{ 1, 2, \ldots, m \},$ the set of ranking
patterns generated by all $\bfx=(x_1,\ldots,x_m)$ satisfying
$x_{i_1}<\cdots<x_{i_m}$ can be obtained from the set of ranking
patterns generated by all $\bfx=(x_1,\ldots,x_m)$ with
$x_{1}<\cdots<x_{m}$ just by relabelling the objects (Lemma
\ref{sigmax} in Section 2). These considerations tell us that it
suffices to consider the case $x_1 < \cdots < x_m.$

Midpoint order depends on the distances between objects. The joint
scale discussed so far is sometimes called the quantitative joint
scale. Another type of joint scale is sometimes considered where
we disregard the metric information of the quantitative joint
scale and take into account only the order of its objects. In this
case we obtain the so-called qualitative joint scale. The set of
admissible rankings of the qualitative joint scale having objects
$x_1,\ldots,x_m$ with $x_{i_1}<\cdots < x_{i_m}$ is, by
definition, the union of the sets of admissible rankings of the
quantitative joint scales whose objects are given by changing only
the distances among $x_1,\ldots,x_m$ while keeping their rank
order $x_{i_1}<\cdots < x_{i_m}.$ Obviously, the number of
qualitative joint scales is $m!/2$ and the number of admissible
rankings of each qualitative joint scale is $2^{m-1}$
(Davison~\cite{dav}). In this paper, we consider quantitative
joint scales exclusively, hence a joint scale always means a
quantitative joint scale.

Suppose the objects $x_1,\ldots,x_m$ are ordered as $x_1< \cdots
<x_m.$ We want to know the number of possible rank orders of the
midpoints $x_{ij}=(x_i+x_j)/2, \ 1 \le i<j \le m.$ Any possible
rank order of the midpoints $x_{ij}, \ 1 \le i<j \le m,$ must
satisfy the condition that the rank $d(i,j)$ of $x_{ij}$  from
left to right on $\bbR$ be increasing in $i$ for any fixed $j$ as
well as increasing in $j$ for any fixed $i.$ Consider the number
$g_m$ of functions $d: \{ (i,j) \mid 1 \le i<j \le m\} \to \{
1,2,\ldots,m(m-1)/2 \}$ satisfying this condition. Clearly $g_m$
serves as an upper bound for the number of possible rank orders of
the midpoints $x_{ij}, \ 1 \le i<j \le m.$ Thrall~\cite{thr}
obtained this number by considering a problem similar to that of
counting the number of standard Young tableaux. However, $g_m$ is
only an upper bound, since the rank order of the midpoints meeting
the above-mentioned condition does not necessarily satisfy other
restrictions induced by the rank order of the objects. Van
Blokland-Vogelesang~\cite{vbl} finds two kinds of such
``intransitive'' midpoint orders by way of ``comparing intervals''
and ``merging intervals.''

In this paper, we find the number of possible rank orders of
midpoints and thereby obtain the number of ranking patterns
generated by the unidimensional unfolding model. This is achieved
by introducing a new type of arrangement called the mid-hyperplane
arrangement. For the general theory of hyperplane arrangements,
see Orlik and Terao~\cite{ort}. Although we give a formula for the
number of ranking patterns for all $m$ in Theorem \ref{1.6}, we
calculate this number only  for $m\le 8$ due to computational
complexity. In addition to determining the number of ranking
patterns, we may ask a further question of interest. Suppose the
$m$ objects are randomly determined. What is the probability that
a given ranking pattern occurs? As will be seen in Section 6, this
problem for $m=5$ reduces to that of finding volumes of some
spherical tetrahedra.

The organization of this paper is as follows.  In Section 2, we
define the mid-hyperplane arrangement and show that the number of
ranking patterns can be obtained by counting the number of
chambers of this arrangement. In Section 3, we reduce the problem
to that of counting the number of points in certain finite sets.
Based on these results, we actually obtain the number of ranking
patterns for $m \le 7$  in Section 4 and for $m=8$ in Section 5.
We also show in those sections that the characteristic polynomial
of the mid-hyperplane arrangement is a product of linear factors
in $\ints[t]$ if and only if $m \le 7.$ In Section 6, we consider
the problem of the probabilities of ranking patterns and give the
answer for $m \le 5$ objects. In Section 7, we mention some open
problems.

\section{Arrangements and ranking maps} Let $m$ be
an integer with $m\geq 3$. Here we define two kinds of hyperplanes in the
$m$-dimensional Euclidean space ${\bbR}^{m}$.

\begin{itemize}
\item[(I)]
$H_{ij} :=\{(x_{1}, \dots, x_{m}) \in {\bbR}^{m} \mid
x_{i}=x_{j}\} \quad (1 \leq i < j \leq m)$.
\end{itemize}
The hyperplane arrangement ${\B}_{m} := \{H_{ij} \mid 1 \leq i < j
\leq m \}$ is called the {\it braid arrangement} \cite[p.13]{ort}.
It has $|{\mathcal B}_{m}| = {\binom{m}{2}}$ hyperplanes. Let
\begin{equation*}
I_{4} :=  \{(p,q,r,s) \mid 1 \leq p < q \leq m, \,\, p < r < s
\leq m,\,\,\, p,q,r,s \,\,\, \text{are distinct} \}.
\end{equation*}

\begin{itemize}
\item[(II)]
$H_{pqrs} :=\{(x_{1}, \dots, x_{m}) \in {\bbR}^{m} \mid
x_{p}+x_{q} = x_{r}+x_{s}\}  \quad (p,q,r,s)\in I_4$.
\end{itemize}
Define the {\it mid-hyperplane arrangement}
\begin{equation*}
{\mathcal A}_{m} := \B_{m} \cup \{H_{pqrs} \mid
(p,q,r,s) \in I_{4} \}.
\end{equation*}
Here $|{\mathcal A}_{m}| = {\binom{m}{2}}+3{\binom{m}{4}}$. For an
arbitrary arrangement $\A$ in $\bbR^{m} $, let
\begin{equation*}
M({\mathcal A}) := {\bbR}^{m} \setminus
\bigcup_{H \in {\mathcal A}}H
\end{equation*}
be the complement of ${\mathcal A}$. The connected components of
$M(\A)$ are called {\it chambers} of  $\A $. Let ${\ch} (\A)$ be
the set of all chambers of $\A$.

Let $\bbP_{m}$ denote the set of all permutations of $\{1, 2,
\dots , m\}$. For $\pi = (i_{1} \dots i_{m}) \in \bbP_{m}$, let
$\hat{\pi}$ denote the corresponding bijection from $\{1, \dots,
m\}$ to itself: $\hat{\pi} (k) = i_{k} \,\,\, (1 \leq k \leq m)$.
In this way we have a one-to-one correspondence between $\bbP
_{m}$ and the symmetric group ${\bbS}_{m},$ which is defined to be
the set of bijections from $\{1, \dots, m\}$ to itself. The group
$\bbS_{m}$ acts on the set $\bbP_{m}$ by
$$\sigma\pi := (\sigma(i_{1}) \dots \sigma(i_{m})) \in \bbP_{m}$$
for $\sigma \in \bbS_{m}$ and $\pi = (i_{1} \dots i_{m}) \in \bbP_{m}$.
The action of $\bbS _{m}$ on $\bbR ^{m}$ is defined by
$$\sigma(x_{1}, \dots, x_{m}) = (x_{\sigma ^{-1} (1)}, \dots,
x_{\sigma ^{-1} (m)}).$$ Then $\bbS_{m} $ acts on $M(\A_{m} )$ and
$M(\B_{m} )$ and therefore on $\ch(\A_{m} )$ and on $\ch(\B_{m}
)$.

It is well known (e.g., Bourbaki~\cite[Ch.5, \S 3, ${\rm
n}^\circ$2, Th.1]{bou}) that the symmetric group $\bbS_{m}$ acts
on $\ch(\B_{m})$ effectively and transitively. In other words, for
any $C, C'\in \ch(\B_{m})$, there exists  a unique $\sigma\in
\bbS_{m} $ with $C' = \sigma C$. In particular, $|\ch(\B_{m})| =
m!$. Let
$$C_{0} := \{(x_{1}, x_{2},\dots , x_{m})
\mid x_{1} < x_{2} <\dots < x_{m} \}$$ be  a  chamber  of the
braid arrangement ${\B}_{m}.$ Then $\ch(\B_{m}) = \{ \sigma C_{0}
\mid \sigma \in \bbS_{m} \}$.

Fix $\bfx = (x_{1}, x_{2} ,\dots , x_{m}) \in M({\mathcal
A}_{m}).$ Plot $m$ points $x_{1}, x_{2} ,\dots , x_{m}$ on the
real line $\bbR$. Let $R(\bfx) := \bbR \setminus \{x_{ij} \mid
1\leq i < j\leq m\}$, where $x_{ij} := (x_{i} + x_{j})/2$ is the
midpoint. Define a map
\[
{\cR_{\bfx}} : R(\bfx) \longrightarrow \bbP_{m}
\]
as follows:
\[
\cR_{\bfx}(y) = (i_{1} i_{2} \dots i_{m}) \Longleftrightarrow
|y - x_{i_{1} } | <
|y - x_{i_{2} } | <
\dots
<|y - x_{i_{m} } |,
\]
where $y \in R(\bfx)$ and $(i_{1} i_{2} \dots i_{m})\in\bbP_{m} $.
The map $\cR_{\bfx} $  is called the {\it ranking map}. The image
of the ranking map $\cR_{\bfx} $ is  the {\it ranking pattern} of
$\bfx\in M(\A_{m} )$.

Suppose $\bfx\in  C_{0} \cap M(\A_{m})$. Then $x_{1} < x_{2}
<\dots < x_{m}$. For $y \in R(\bfx)$ and $1 \leq i < j \leq m$, we
have
$$y < x_{ij} \Longleftrightarrow |y-x_{i}|<|y-x_{j}|
\Longleftrightarrow i \text{~precedes~} j \text{~in~}
\cR_{\bfx}(y),$$
$$y > x_{ij} \Longleftrightarrow |y-x_{i}|>|y-x_{j}|
\Longleftrightarrow j \text{~precedes~} i \text{~in~}
\cR_{\bfx}(y).$$ Imagine that the point $y$ moves on the real line
$\bbR$ from left to right. When $y$ is sufficiently small,
$\cR_{\bfx}(y) = (12 \dots m)$. Every time  $y$ ``passes''
$x_{ij}$, the two integers $i$ and $j$, which are adjacent  in
$\cR_{\bfx}(y)$, switch their positions. When $y$ is sufficiently
large, $\cR_{\bfx}(y) = (m \dots 21)$.

\begin{example}
Let $m = 3$
and $x_{1} < x_{2} < x_{3}$.
  Then
$$
\cR_{\bfx} (y) =
\begin{cases}

(123) \,\, \text{if~} y < x_{12},\\
(213) \,\, \text{if~}  x_{12} < y < x_{13} ,\\
(231) \,\, \text{if~}  x_{13} < y < x_{23} ,\\
(321) \,\, \text{if~}  x_{23} < y .
\end{cases}
$$
\end{example}

\begin{lemma}
\label{sigmax} Let $\sigma \in \bbS _{m}$, $\bfx \in M(\A_{m})$
and $y \in R(\bfx)$. Then $$\cR_{\sigma \bfx}(y) = \sigma
(\cR_{\bfx}(y)).$$
\end{lemma}

\begin{proof}
Suppose $\bfx = (x_{1}, \dots, x_{m})$.
Then
\begin{multline*}
\cR_{\sigma \bfx}(y) =(i_{1} \dots i_{m})
\Longleftrightarrow |y-x_{\sigma^{-1}(i_{1})} | < \dots <
|y-x_{\sigma^{-1}(i_{m})} |\\
\Longleftrightarrow \cR_{\bfx}(y) = (\sigma^{-1} (i_{1}) \dots
\sigma^{-1} (i_{m})) \Longleftrightarrow \sigma(\cR_{\bfx}(y)) =
(i_{1} \dots i_{m}).
\end{multline*}
\end{proof}

\begin{lemma}
\label{1.5}
Let $\sigma\in \bbS_{m} $ and $\bfx, \bfx'\in \sigma C_{0} \cap M(\A_{m} )$.
Then
$\bfx$ and $\bfx'$ lie in the same chamber of $\A_{m} $
if and only if
the following statement holds true:
$$
x_{pq} > x_{rs} \Longleftrightarrow x'_{pq} > x'_{rs}
$$
for each $(p,q,r,s)\in I_{4}.$
\end{lemma}

\begin{proof}
Each chamber of $\A_{m}$ inside $\sigma C_{0} $ is equal to the
intersection of $\sigma C_{0} $ and half-spaces defined by either
$2(x_{pq} - x_{rs} ) = x_{p} + x_{q} - x_{r} - x_{s} > 0$ or
$2(x_{pq} - x_{rs} ) = x_{p} + x_{q} - x_{r} - x_{s}  < 0$ for
$(p,q,r,s)\in I_{4}$.
\end{proof}

\begin{theorem}
\label{1.3} Let $\sigma\in \bbS_{m} $ and $\bfx, \bfx'\in \sigma
C_{0} \cap M(\A_{m} )$.  Then $\bfx$ and $\bfx'$ have the same
ranking pattern if and only if $\bfx$ and $\bfx'$ lie in the same
chamber of $\A_{m}. $
\end{theorem}

\noindent
\begin{proof}
Assume first that $\sigma= 1$, so $\bfx, \bfx' \in C_{0}\cap
M(\A_{m} ).$

Suppose that $\bfx$ and $\bfx'$ lie in the same chamber of $\A_{m}
$. Write $\bfx' = (x'_{1}, x'_{2}, \dots , x'_{m})$ and $x'_{ij}
:= (x'_{i} + x'_{j} )/2\,\, (1\leq i<j\leq m)$. By Lemma
\ref{1.5}, we have
\[
x_{i_{1} j_{1} }
<
x_{i_{2} j_{2} }
<
\dots
<
x_{i_{t} j_{t} }, \,\,\,\,
x'_{i_{1} j_{1} }
<
x'_{i_{2} j_{2} }
<
\dots
<
x'_{i_{t} j_{t} },
\]
where $t = \binom{m}{2}$.
This shows
\[
\im \cR_{\bfx} = \{ \pi_{0}, \pi_{1}, \dots, \pi_{t} \}=
\im \cR_{\bfx'},
\]
where $\pi_{0}, \pi_{1}, \dots \pi_{t} \in \bbP _{m}$ are defined
inductively by
\begin{eqnarray*}
\pi_{0} &=& (12\dots m),\\
\pi_{s} &=& [i_{s}j_{s}]\pi _{s-1} \,\,\,(1 \leq s \leq t).
\end{eqnarray*}
Here $[ij] \in \bbS_{m}\,\,
(1 \leq i < j \leq m)$
denotes the transposition of $i$ and $j$.

Conversely, assume $ \im\cR_{\bfx} = \im\cR_{\bfx'}. $ For
$\pi=(i_{1} i_{2} \dots i_{m} )\in \bbP_{m} $, let $\iota(\pi)$
denote the number of inversions in $\pi$:
\[
\iota(\pi) := |\{(k_{1}, k_{2}  ) \mid
k_{1} < k_{2}, i_{k_{1} } > i_{k_{2} }  \}|.
\]
  As the point $y$ moves on the real line from left to right,
  $\iota(\cR_{\bfx}(y))$ increases one by one. So we may write
  \[
  \im \cR_{\bfx} = \im \cR_{\bfx'} = \{\pi_{0}, \pi_{1}, \dots , \pi_{t}\}
  \]
  such that $\iota(\pi_{s} ) = s,\,\, (0\leq s \leq t)$.
Also there exists a unique transposition $[i_{s}j_{s}]$ such that
$\pi_{s} = [i_{s} j_{s}] \pi_{s-1} \,\, (1\leq s \leq t).$ Thus $
x_{i_{1}j_{1}} < x_{i_{2}j_{2}} < \dots < x_{i_{t}j_{t}}$ and $
x'_{i_{1}j_{1}} < x'_{i_{2}j_{2}} < \dots < x'_{i_{t}j_{t}}. $ It
follows from Lemma \ref{1.5} that $\bfx$ and $\bfx'$ lie in the
same chamber of $\A_{m} $.

For a general $\sigma \in \bbS_{m}$, let $\bfy
:=\sigma^{-1}\bfx\in C_{0} \cap M(\A_{m})$ and $\bfy' :=
\sigma^{-1}\bfx'\in C_{0} \cap M(\A_{m})$. By Lemma \ref{sigmax},
\begin{multline*}
\im \cR_{\bfx}
=
\im \cR_{\bfx'}
\Leftrightarrow
\sigma^{-1} (\im \cR_{\bfx})
=
\sigma^{-1} (\im \cR_{\bfx'})
\Leftrightarrow
\im \cR_{\sigma^{-1} \bfx}
=
\im \cR_{\sigma^{-1} \bfx'} \\
\Leftrightarrow
\im \cR_{\bfy}
=
\im \cR_{\bfy'}
\Leftrightarrow
\bfy \text{~and~}\bfy'
\text{~lie in the same chamber of~}
\A_{m} \\
\Leftrightarrow
\bfx \text{~and~}\bfx'
\text{~lie in the same chamber of~}
\A_{m}.
\end{multline*}
\end{proof}

Let $r(m)$ denote the number of
ranking patterns when
$\bfx$ runs over the set $C_{0} \cap M(\A_{m} )$:
\[
r(m) := |\{\im\cR_{\bfx}\mid \bfx\in C_{0}\cap M(\A_{m} )  \}|.
\]
Note that for each $\sigma \in \bbS_m,$
$|\{\im\cR_{\bfx}\mid \bfx \in \sigma C_{0}\cap M(\A_{m} )  \}|$
is equal to $r(m)$ by Lemma \ref{sigmax}.

\begin{theorem}
\label{1.6}
$r(m) = |\ch(\A_{m} )|/(m!)$.
\end{theorem}

\begin{proof}
By Theorem \ref{1.3}, $r(m)$ is equal to the number of chambers of
$\A_{m} $ which lie inside $C_{0} $. Thus we have $ |\ch(\A_{m} )|
= r(m) |\bbS_{m} | = r(m) (m!).$
\end{proof}

\section{The number of chambers of $\A_{m}$ }

In this section, we study the number $|\ch(\A_{m})|$ of chambers
of $\A_{m} $.

First let us review some general results about the number of
chambers and the characteristic polynomial. Let $\bbK$ be a field
and $V$ an $\ell$-dimensional vector space over $\bbK$. Assume
that $\A$ is an arbitrary arrangement of hyperplanes in $V$.  Let
$L=L(\A)$ be the set of nonempty intersections of elements of
$\A$. An element  $X \in L$ is called an {\it edge} of $\A$.
Define a {\it partial order} on $L$ by $ X \leq Y
\Longleftrightarrow Y \subseteq X$. Note that this is reverse
inclusion. Thus $V$ is the unique minimal element of $L$.

Let $\mu: L  \ra \ints$ be the M\"obius function
of $L$ defined by
$\mu(V)=1$, and for $X>V$ by the recursion
$$\sum_{Y \leq X}\mu(Y)=0.$$
The {\it characteristic polynomial} of $\A$ is
$$\chi(\A,t)=\sum_{X \in L}\mu(X)t^{\dim X}.$$
The next two theorems give geometric meaning to special values of
the characteristic polynomial.

 \begin{theorem}[Zaslavsky \cite{zas}]
 \label{2.1}
 If $\bbK=\bbR$, then $|\chi(\A, -1)|=|\ch(\A)|.$
 \end{theorem}

An arrangement $\A$ is called {\it essential} if the dimension of
a maximal element of $L(\A)$ is zero. The mid-hyperplane
arrangement $\A_{m}$ is not essential because the line $l= {\rm
span}\{ {\bf 1} \} = \{ \lambda{\bf 1} \mid \lambda \in \bbR \}
\subset \bbR^m,$ where ${\bf 1} \in \bbR^m$ is the vector of 1's
is a maximal element. This implies that $\chi(\A_m,t)$ is
divisible by $t$. The fact that $l$ is contained in every
hyperplane of $\A_m$ implies that $\chi(\A_m,t)$ is also divisible
by $(t-1)$. Thus $\chi(\A_m,t)/t(t-1)$ is a monic polynomial of
degree $m-2$.

Let $H_{0} $ be the hyperplane defined by $x_{1} = 0$. Define
$\A_{m}^{*}  := \A_{m}  \cup \{H_{0} \}$. Then $\A_{m}^{*} $ is
essential and the lattice $L(\A_{m})$ is isomorphic to the
sublattice defined by $ L(\A_{m}^{*})_{\geq H_{0} } := \{X \in
L(\A_{m}^{*}) \mid X \geq H_{0} \}$.

 \begin{theorem}
[Crapo-Rota \cite{crr}, Terao \cite{tertohoku} (4.10)] \label{2.2}
Let $\bbF_{q}$ be a finite field of $q$ elements. If
$\bbK=\bbF_{q}$, then $\chi(\A, q) =|M(\A)|.$
 \end{theorem}

When $\bbK=\bbF_{q}$ and $V$ is a finite set of $q^\ell$ elements,
$\chi(\A, q)$ can be evaluated by counting the number of points
not on any hyperplane $H\in {\cal A}$ in $V$. Let $q$ be a prime
number greater than $m$. Let $\A_{m, q}^{*}$ be the modulo $q$
reduction of $\A_{m}^{*}$ in $(\mathbb Z_{q} )^{m}$. In other
words, the hyperplanes belonging to $\A_{m, q}^{*}$ are:

\medskip

$(0)$ $H_{0} :=\{(x_{1}, \dots, x_{m}) \in (\mathbb Z_{q} )^{m}
\mid x_{1} = 0 \} $,

\smallskip

$(I_{q})$
$H_{ij} :=\{(x_{1}, \dots, x_{m}) \in (\mathbb Z_{q} )^{m} \mid
x_{i}=x_{j}\} \quad (1 \leq i < j \leq m)$,
and
\smallskip

$(II_{q})$
$H_{pqrs} :=\{(x_{1}, \dots, x_{m}) \in (\mathbb Z_{q} )^{m} \mid
x_{p}+x_{q} = x_{r}+x_{s}\}  \quad
(p,q,r,s\in I_{4}).$

\medskip

The arrangement $\A_{m, q}^{*}$ is essential. The modulo $q$
reduction $\A_{m, q} $ of $\A_{m} $ is composed of the hyperplanes
of type $(I_{q})$ and $(II_{q} )$ above. Note that the lattice
$L(\A_{m, q})$ is isomorphic to the sublattice defined by $
L(\A_{m, q}^{*})_{\geq H_{0} } := \{X \in L(\A_{m, q}^{*} ) \mid X
\geq H_{0} \}$. Therefore, the intersection lattices $L(\A_{m})$
and $L(\A_{m, q} )$ are isomorphic if  $L(\A_{m}^{*})_{\geq H_{0}
}$ and $L(\A_{m, q}^{*})_{\geq H_{0} }$ are isomorphic.

Let $C$ be the coefficient matrix of $\A_{m}^{*} $. For example
when  $m = 4$,
$$C
=
\left(
\begin{array}{cccccccccc}
1&1&1&1&0&0&0&1&1&1\\
0&-1&0&0&1&1&0&-1&-1&1\\
0&0&-1&0&-1&0&1&-1&1&-1\\
0&0&0&-1&0&-1&-1&1&-1&-1
\end{array}
\right).
$$

Consider the $m$-minors of $C$. Each $m$-minor is parametrized by
the set of $m$ columns used for the minor. It is known that the
intersection lattice of an essential arrangement is completely
determined by the information which $m$-minors vanish and which do
not \cite[Proposition 3]{termoduli}. Thus we have

\begin{theorem}
\label{2.3} Define $f(m) := \max \{ |\det T| \mid T $ is an
$m$-minor of~$C$ and $T$ contains the column $(1, 0, \dots
,0)^{T}$ as its first column$ \}.$ Let $q$ be a prime number
greater than $f(m)$. Then $L(\A_{m})$ and $L(\A_{m, q})$ are
isomorphic.
\end{theorem}

Next we will find an upper bound for $f(m)$. Let $m\geq 3$ as
always. Consider the following three conditions concerning a
matrix:

(i) every entry of the matrix is either $-1$, $0$ or $1$,

(ii) in every column $1$ appears at most twice,

(iii) in every column $-1$ appears at most twice.

\noindent Define $$g(m) := \max\{|\det A| \mid A \text{~is an~}
(m-1)\times (m-1) \text{-matrix satisfying (i,ii,iii)}\}.$$  It is
clear
 that $f(m) \leq g(m)$.

\begin{lemma}
\label{2.4}
$$
g(m) = 2^{m-2} \mbox{ for } m\leq 5, \quad g(m) \leq 8 \cdot
3^{m-5} \mbox{ for }m\geq 6.
$$
\end{lemma}

\begin{proof}
We argue by induction on $m$. For $m = 3, 4$ direct computation
shows the result. The values are attained by
\[
g(3)
=
2=
\det
\left(
\begin{array}{cc}
1&-1\\
1&1
\end{array}
\right),
\,\,\,
g(4)
=
4
=
\det
\left(
\begin{array}{ccc}
-1&1&-1\\
-1&-1&1\\
 0&1&1
\end{array}
\right).
\]
Suppose $m=5$. We show that $g(5) = 8$. Note that
\[
\det
\left(
\begin{array}{cccc}
1&1&1&1\\
1&-1&-1&1\\
0&-1&1&-1\\
0&1&-1&-1
\end{array}
\right)
=
8.
\]
Thus $g(5)\geq 8$. We must show $ g(5)=\det A\leq 8$. Let $\mathbf
a_{i} $ be the $i$th column of $A$. Denote the number of nonzero
elements in $\mathbf a_{i} $ by $val(\mathbf a_{i} )$. If $A$ has
a column $\mathbf a_{i}$ with $val(\mathbf a_{i})\leq 2$, then
$g(5) \leq 2 g(4) = 8$. So we may assume that $3 \leq val(\mathbf
a_{i})\leq 4$ for every $i$.   Define $ \mathbf a^{*}_{i} $ to be
the uniquely determined four-dimensional column vector with two
$1$'s and two $-1$'s which satisfies the following property:
\begin{eqnarray*}
\mbox{if~} val(\mathbf a_{i})&=&4, \mbox{~then~} \mathbf a^{*}_{i}
=
\mathbf a_{i}  , \\
\mbox{if~} val(\mathbf a_{i})&=&3, \mbox{~then~} \mathbf a^{*}_{i}
\mbox{~is obtained from~} \mathbf a_{i}\mbox{~by replacing } \\ &&
\shoveleft{\mbox{~the unique zero in~} \mathbf a_{i}\mbox{ by
either $1$ or $-1$}.}
\end{eqnarray*}
For example,
\[
\mbox{~if~}
\mathbf a_{i}
=
\left(
\begin{array}{c}
1\\
-1\\
0\\
-1
\end{array}
\right) \mbox{ then }
\mathbf a^{*}_{i} = \left(
\begin{array}{c}
1\\
-1\\
1\\
-1
\end{array}
\right).
\]
Since $\binom{4}{2}/2 = 3 < 4$,  among the four vectors $\mathbf
a^{*}_{i} \,\,(1\leq i\leq 4)$ at least two are either equal to or
the negative of each other. Without loss of generality, we may
assume $\mathbf a^{*}_{1} =  \mathbf a^{*}_{2}$. Then $\mathbf
a_{12} := \mathbf a_{1}-\mathbf a_{2}$ is composed only of $0, -1$
and $1$ with $val(\mathbf a_{12} )=2$. We get
\begin{multline*}
g(5)=
\det A =
\det
(\mathbf a_{1}, \mathbf a_{2}, \mathbf a_{3}, \mathbf a_{4})
=
\det
(\mathbf a_{1}-\mathbf a_{2}, \mathbf a_{2}, \mathbf a_{3}, \mathbf a_{4})
\\
=
\det
(\mathbf a_{12}, \mathbf a_{2}, \mathbf a_{3}, \mathbf a_{4})
\leq
val(\mathbf a_{12} ) g(4) = 2 g(4) = 8.
\end{multline*}

For the induction step, assume $m\geq 6$. Choose an $(m-1)\times
(m-1)$-matrix $A$ satisfying (i), (ii) and (iii) with $\det A =
g(m)$.  Use the Laplace expansion formula along $\mathbf a_{i} $
to get
\[
g(m) \leq val(\mathbf a_{i} ) g(m-1)
\]
for each $i$. If $val(\mathbf a_{i} ) = 4$ for every column
$\mathbf a_{i} $ of $A$, then all the rows sum up to the zero
vector and thus $\det A = 0$. This is a contradiction. Therefore
we may assume that there exists a column $\mathbf a_{i} $ with
$val(\mathbf a_{i} )\leq 3$; so we obtain $g(m) \leq 3 g(m-1)$.
\end{proof}
Theorem \ref{2.3} and Lemma \ref{2.4} imply
\begin{theorem}
\label{2.5}
If a prime number $q$ satisfies
\[
q >
\begin{cases}
2^{m-2}    \mbox{~~~~~~~if } m\leq 5,\\
8\cdot 3^{m-5} \mbox{~~~~if } m\geq 6,
\end{cases}
\]
then
the intersection lattices
$L(\A_{m} )$ and $L(\A_{m, q})$ are isomorphic
and
\[
\chi(\A_{m},  q) = |M(\A_{m, q})|.
\]
\end{theorem}

The following theorem shows that we can fix $x_1=0$, $x_2=1$  in
counting $|M(\A_{m, q})|$.

\begin{theorem}
\label{2.7}
Define
\[
M_{1} (m, q) := \{(0, 1, x_{3}, \dots , x_{m}  )\in M(\A_{m, q})\}.
\]
Under the assumption of Theorem \ref{2.5}, we have
\[
\frac{\chi(\A_{m}, q)}{q(q-1)} = |M_{1} (m, q)|.
\]
 \end{theorem}

\begin{proof}
Consider the action of the additive group $\bbF_{q}$ on $M(\A_{m,
q})$ by
\[
(x_{1}, x_{2}, \dots , x_{m}) \mapsto
(x_{1} + \alpha, x_{2} + \alpha, \dots , x_{m} + \alpha)
\,\,\, (\alpha\in\bbF_{q} ).
\]
The set of orbits under this action is represented by
the set
$$M_{0} := \{(0, x_{2} , x_{3}, \dots , x_{m}  )
\in M(\A_{m, q})\}.
$$
Thus $|M_{0} | = \chi(\A_{m}, q )/q$. Next consider the action of
the multiplicative group $\bbF_{q}^{\times} := \bbF \setminus
\{0\} $ on $M_{0} $ by
\[
(0, x_{2}, \dots , x_{m}) \mapsto
(0, x_{2}  \beta, \dots , x_{m}  \beta)
\,\,\, (\beta\in\bbF_{q}^{\times}).
\]
The set of orbits under this action is represented by
the set $M_{1} (m, q)$.
Thus
$|M_{1} (m, q)| = |M_{0} |/(q-1) = \chi(\A_{m}, q )/q(q-1).$
\end{proof}

Our method to count the number of chambers of $\A_{m} $ is as
follows: let $q_{i}$ $(i=1,\dots,m-2)$ be primes satisfying the
conditions of Theorem \ref{2.5}. Count the number of points in the
set $M_{1} (m, q_{i} )$ for each $i$. By Theorem \ref{2.7}, we
have $\chi(\A_{m}, q_{i})/q_{i} (q_{i} - 1) = |M_{1} (m, q_{i}
)|$. Since $\chi(\A_{m}, t)/t (t - 1)$ is a monic polynomial of
degree $m-2$, the data $|M_{1} (m, q_{i} )|\,\,(i=1,\dots, m-2)$
determine the characteristic polynomial $\chi(\A_{m}, t)$ of
$\A_{m} $. Theorem \ref{2.1} asserts
\[
|\ch(\A_{m} )| = |\chi(\A_{m}, -1)|.
\]
The number $r(m)$ of the ranking patterns when $\bfx$ runs over
the set $C_{0} \cap M(\A_{m})$ is obtained by $r(m) =
|\chi(\A_{m}, -1)|/(m!)$ by Theorem \ref{1.6}.

\section{The number of ranking patterns for $m \leq 7$}

In this section we determine $\chi(\A_{m}, t),$ $|\ch(\A_{m})|$
and $r(m)$ for $m\leq 7$. The case $m = 3$ is known because
$\A_{3} = \B_{3}$. Let $m=4$.  If  $q >4$ is a prime, then Theorem
\ref{2.7} gives
\[
\frac{\chi(\A_{4}, q)}{q(q-1)} = |M_{1} (4, q)|.
\]
Let $$p(t) := \frac{\chi(\A_{4}, t)}{t(t-1)}.$$ Then $p(t)$  is a
monic quadratic polynomial. We find
     \begin{equation*}
     p(5) = |M_{1} (4, 5)|=0 \,\,\, \text{and} \,\,\,   p(7) = |M_{1} (4, 7)|=8.
     \end{equation*}
Theorems \ref{2.1} and \ref{1.6} give
\begin{gather*}
      p(t) = t^{2} - 8t + 15 = (t-3)(t-5), \,\,
\chi(\A_{4}, t ) = t(t-1)(t-3)(t-5),\\
|\ch(\A_{4})| = 48, \,\, \text{and} \,\,\, r(4) = 2.
\end{gather*}
Using the same method, computer calculations provide the following
table:

\begin{theorem}
\label{3.1}
\[
\begin{tabular}{c|c|c|c}
$m$ & $\chi(\A_{m}, t )$ & $|\ch(\A_{m})|$ & $r(m)$\\
\hline
$3$ & $t(t-1)(t-2)$ & $6$ & $1$ \\
$4$ & $t(t-1)(t-3)(t-5)$ & $48$ & $2$ \\
$5$ & $t(t-1)(t-7)(t-8)(t-9)$ & $1440$ & $12$ \\
$6$ & $t(t-1)(t-13)(t-14)(t-15)(t-17)$ & $120960$ & $168$ \\
$7$ & $t(t-1)(t-23)(t-24)(t-25)(t-26)(t-27)$ & $23587200$ & $4680$
\end{tabular}
\]
\end{theorem}

 \begin{corollary}
   \label{3.2}
 If $m \leq 7$, then the characteristic polynomial
$\chi(\A_{m}, t )$ is a product of linear factors in $\ints[t]$.
   \end{corollary}

     \medskip
     \noindent
{\em Remark.} Define $a_n=n(n^{n-1}-1)((n-2)!)/(n-1)$. We note
that $r(m)=a_{m-2}$ for $m = 3, 4, 5, 6, 7$ but we do not have any
reasonable interpretation for the coincidence at this writing.

\section{The number of ranking patterns for $m \geq 8$}

We determine $r(8)$ first. Evaluating $r(m)$ for $m \ge 9$ is not
feasible  at present with our brute force counting method. Next we
prove a theorem about the characteristic polynomial $\chi(\A_{m},
t )$ for $m\geq 8$.

For $m=8$ we used a computer to count $|M_{1} (8, q)| $ with the
primes $q=$ 223, 227, 229, 233, 239, 241, all  greater than
$8\cdot 3^{8-5} = 216$. Theorem \ref{2.7} implies:
\begin{theorem}
\label{4.2}
\begin{gather*}
\chi(\A_{8}, t) = t(t-1)
(t-35)(t-37)(t-39)(t-41)(t^2-85t+1926),\\
|\ch(\A_{8})| = 9248117760,\\
r(8) = 229386.
\end{gather*}
\end{theorem}

\medskip
\noindent {\em Remark.}  The coincidence of $r(m)$ and $a_{m-2}$
does not hold for $m = 8$. Here $r(8) = 229386 > a_{6} = 223920.$

Write
\[
\chi(\A_{m}, t ) = \sum_{k=0}^{m} \mu_{k} t^{m-k}.
\]
It is known that
$$\mu_{0}=1, \,\, \mu_{1} = -|\A_{m}| = -\binom{m}{2} - 3 \binom{m}{4},
\,\,   \mu_{m} = 0.$$
Although we do not have a general formula for $\mu_{k} $,
routine calculations yield a formula for $\mu_{2} $:
\begin{theorem}
\label{4.1}
$$
\mu_{2}
=2\binom{m}{3}+15\binom{m}{4}+120\binom{m}{5}
+375\binom{m}{6}+630\binom{m}{7}+315\binom{m}{8}.
$$
\end{theorem}

\begin{theorem}
\label{4.3} The characteristic polynomial $ \chi(\A_{m}, t)$ is a
product of linear factors in $\ints[t]$ if and only if $m\leq 7$.
\end{theorem}

\begin{proof}
This follows from Corollary \ref{3.2} when $m\leq 7$. Let $m\geq
8$. Suppose that the characteristic polynomial is a product of
linear factors in $\ints[t]$:
\[
\chi(\A_{m}, t ) = \sum_{k=0}^{m} \mu_{k} t^{m-k}
= t (t-1) (t-b_{2})\dots (t-b_{m-1})
\]
 for $b_{2}, \dots , b_{m-1} \in \ints$.

Applying Theorem \ref{4.1}, we have
\begin{align*}
\sum_{i=2}^{m-1}b_i &= -\mu_{1} -1
=|\A_{m} | - 1= -1+\binom{m}{2}+3\binom{m}{4}, \\
\sum_{2\le i<j\le m-1}b_ib_j
&=\mu_2-\sum_{i=2}^{m-1}b_i \\
&=1-\binom{m}{2}+2\binom{m}{3}+12\binom{m}{4} \\
& \ \ \ \ +120\binom{m}{5}
+375\binom{m}{6}+630\binom{m}{7}+315\binom{m}{8}.
\end{align*}

Therefore
\begin{align*}
\sum_{i=2}^{m-1}\left( b_i-\frac{\sum_{i=2}^{m-1}b_i}{m-2}\right)^2
&=\sum_{i=2}^{m-1}b_i^2-\frac{\left( \sum_{i=2}^{m-1}b_i\right)^2}{m-2} \\
&=\left( \sum_{i=2}^{m-1}b_i \right)^2 -2\sum_{2\le i<j\le m-1}b_ib_j
-\frac{\left( \sum_{i=2}^{m-1}b_i \right)^2}{m-2} \\
&=\frac{(m-3)\left( \sum_{i=2}^{m-1}b_i \right)^2}{m-2}
-2\sum_{2\le i<j\le m-1}b_ib_j.
\end{align*}

Compute
\begin{align*}
h(m)
&:=(m-2)\sum_{i=2}^{m-1}\left(b_i-\frac{\sum_{i=2}^{m-1}b_i}{m-2}\right)^2\\
&=(m-3)\left\{ -1+\binom{m}{2}+3\binom{m}{4}\right\}^2
-2(m-2)\left\{ 1-\binom{m}{2}+2\binom{m}{3}\right.\\
& \ \ \ \ \left.+12\binom{m}{4}+120\binom{m}{5}+375\binom{m}{6}
+630\binom{m}{7}+315\binom{m}{8} \right\} \\
&=1+\frac{98m}{3}-\frac{1573m^2}{16}+\frac{5423m^3}{48}
-\frac{12787m^4}{192}+\frac{527m^5}{24}-\frac{391m^6}{96}\\
& \ \ \ \ \ +\frac{19m^7}{48}-\frac{m^8}{64}.
\end{align*}
Thus $h(m) \geq 0$ for $m > 2$. On the other hand, we may check by
standard calculus techniques that $h(m) < 0$ whenever $m \geq 8$.
This is a contradiction.
\end{proof}

\section{Probabilities of ranking patterns}
We counted the number of possible ranking patterns in the
preceding sections. Here we investigate the probabilities of
ranking patterns when the objects $x_1,\ldots, x_m$ are randomly
determined. For $m=4,$ the problem is trivial by symmetry
considerations as long as the four objects are independently and
identically distributed.

We consider the case $m=5$ and assume that $\bfx=(x_1,\ldots,x_5)
\in \bbR^5$ is distributed according to an arbitrary spherical
distribution. Note that $\bfx \in M({\mathcal A}_5)$ with
probability one. For $m=5,$ there are $1440$ possible ranking
patterns in all. By relabelling the indices it suffices to
consider the case
\begin{equation}
\label{restriction1}
x_1< \cdots <x_5.
\end{equation}
Furthermore, by replacing $x_{i} $ by $-x_{i} $, it suffices to
consider the case
\begin{equation}
\label{restriction2}
x_1< \cdots <x_5, \ \ x_{24} < x_{15}.
\end{equation}

Under restriction (\ref{restriction2}), we have $1440/(5! \cdot
2)=r(5)/2=6$ possible ranking patterns, which are characterized by
the following midpoint orders (Lemma \ref{1.5}, Theorem
\ref{1.3}):
\begin{eqnarray*}
&({\rm I})&   x_{14}<x_{23}<x_{24}<x_{15}<x_{25}<x_{34},
\label{mporder(I)} \\
&({\rm II})&  x_{14}<x_{23}<x_{24}<x_{15}<x_{34}<x_{25},  \\
&({\rm III})& x_{14}<x_{23}<x_{24}<x_{34}<x_{15}<x_{25}, \\
&({\rm IV})& x_{23}<x_{14}<x_{24}<x_{15}<x_{25}<x_{34}, \\
&({\rm V})&  x_{23}<x_{14}<x_{24}<x_{15}<x_{34}<x_{25}, \\
&({\rm VI})& x_{23}<x_{14}<x_{24}<x_{34}<x_{15}<x_{25}.
\label{mporder(VI)}
\end{eqnarray*}

We are interested in the conditional probabilities of the six
midpoint orders above assuming (\ref{restriction2}). Recall that
these midpoint orders represent chambers of ${\mathcal A}_5$
(Lemma \ref{1.5}). We argue next that our problem reduces to
computing the spherical volumes of the restrictions of some
chambers of ${\mathcal A}_5$ to the three-dimensional unit sphere.

We begin by recalling that all hyperplanes in ${\mathcal A}_5$
contain the line $l= {\rm span}\{ {\bf 1} \} = \{ \lambda{\bf 1}
\mid \lambda \in \bbR \} \subset \bbR^5,$ where ${\bf 1} \in
\bbR^5$ is the vector of 1's. The orthogonal projection of
$\bfx=(x_1,\ldots,x_5) \in \bbR^5$ onto $H_0'= l^{\perp} =\{ (x_1,
\ldots, x_5) \in \bbR^5 \mid x_1+\cdots+x_5=0 \}$ will be denoted
by $\bfz:=(x_1-\bar{x}, \ldots,x_5-\bar{x}),$ where
$\bar{x}=(x_1+\cdots+x_5)/5.$ Since $\bfx$ is assumed to be
distributed as a spherical distribution, the marginal distribution
of the orthogonal projection $\bfz$ is a spherical distribution of
one less dimension (Muirhead~\cite[p.34]{mui}). Now, any $\bfx \in
M({\mathcal A}_5)$ and its orthogonal projection $\bfz$ are on the
same side of each hyperplane in ${\mathcal A}_5,$ so for any
chamber $ C  \in {\rm \bf Ch}({\mathcal A}_5),$ we have ${\rm
Prob}(\bfx \in  C) ={\rm Prob}(\bfz \in  C_{H_0'})$ with $
C_{H_0'}:= C\cap H_0'.$ This $ C_{H_0'}$ can be regarded as a
chamber of the arrangement ${\mathcal A}_5' :=\{ H \cap H_0' \mid
H \in {\mathcal A}_5 \}$ in $H_0'.$

Each hyperplane in ${\mathcal A}_5'$ contains the origin. Thus its
chambers are the interiors of polyhedral cones in $H_0'.$ As a
result, for each $ C_{H_0'} \in {\rm \bf Ch}({\mathcal A}_5'),$ we
have that $\bfz \in  C_{H_0'}$ is equivalent to $\bfz/\| \bfz \|
\in  C_{\bbS^3}:= C_{H_0'} \cap \bbS^3,$ where  $\bbS^3:= \{ (x_1,
\ldots, x_5) \in H_0' \mid x_1^2+\cdots +x_5^2=1 \}$ is the unit
sphere in $H_0'.$ Together with the uniformity of the distribution
of $\bfz/\| \bfz \|$ on $\bbS^3,$ this yields $ {\rm Prob}(\bfz
\in  C_{H_0'}) ={\rm Prob}(\bfz/\| \bfz \| \in  C_{\bbS^3}) ={\rm
Vol}( C_{\bbS^3})/{\rm Vol}(\bbS^3). $

We conclude that for any chamber $ C$ of ${\mathcal A}_5,$
\[
{\rm Prob}\left(\bfx \in  C\right)
=\frac{{\rm Vol}( C_{\bbS^3})}{{\rm Vol}(\bbS^3)}
\]
with $ C_{\bbS^3}= C \cap \bbS^3.$ Thus the probability of $\bfx$
being in chamber $ C \in {\rm \bf Ch}({\mathcal A}_5)$ is
proportional to the volume of $ C_{\bbS^3}= C\cap \bbS^3.$
Therefore, the desired conditional probabilities under
(\ref{restriction2}) are given by the ratios of the volumes of the
chambers $ C_{\bbS^3}$ corresponding to the six midpoint orders to
the volume of the union $T:= \{ (x_1, \ldots, x_5) \mid x_1 \le
\cdots \le x_5, \ x_{24} \le x_{15} \} \cap \bbS^3$ of their
closures.

The binding inequalities of the spherical chambers associated with
the six midpoint orders are
\begin{eqnarray*}
&({\rm I})& x_{14}<x_{23}, \ x_{25}<x_{34}, \ x_3<x_4, \ x_{24}<x_{15}, \\
&({\rm II})& x_{15}<x_{34}, \ x_{14}<x_{23}, \ x_{24}<x_{15}, \ x_3<x_4, \ x_{34}<x_{25}, \\
&({\rm III})& x_{14}<x_{23}, \ x_2<x_3, \ x_3<x_4, \ x_{34}<x_{15}, \\
&({\rm IV})& x_1<x_2, \ x_{25}<x_{34}, \ x_{23}<x_{14}, \ x_{24}<x_{15}, \\
&({\rm V})& x_{15}<x_{34}, \ x_{23}<x_{14}, \ x_{24}<x_{15}, \ x_{34}<x_{25}, \\
&({\rm VI})& x_1<x_2, \ x_2<x_3, \ x_{23}<x_{14}, \ x_{34}<x_{15}.
\end{eqnarray*}
With the exception of (II), the closures of these chambers are
spherical tetrahedra
\begin{eqnarray*}
&({\rm I})& FBGH,  \\
&({\rm III})&  AFED, \\
&({\rm IV})&  FBGC, \\
&({\rm V})&  CGFE, \\
&({\rm VI})&  AFCE
\end{eqnarray*}
where
\begin{eqnarray*}
&& A=(-1,-1,-1,-1,4)/\sqrt{20},
\ \ B=(-3,-3,2,2,2)/\sqrt{30}, \\
&& C=(-2,-2,-2,3,3)/\sqrt{30},
\ \ D=(-1,0,0,0,1)/\sqrt{2}, \\
&& E=(-7,-2,-2,3,8)/\sqrt{130},
\ \ F=(-4,-4,1,1,6)/\sqrt{70}, \\
&& G=(-2,-1,0,1,2)/\sqrt{10},
\ \ H=(-8,-3,2,2,7)/\sqrt{130};
\end{eqnarray*}
Chamber (II) is a quadrilateral pyramid $FEDHG,$ which can be
divided into two tetrahedra, say, $FEDG$ and $FDGH.$ Note that
this observation implies that the closures of the chambers of the
mid-hyperplane arrangement ${\mathcal A}_m$ are not necessarily
simplices. See Figures 1 and 2.

The volumes of the seven
spherical tetrahedra mentioned above
can be computed as
\begin{eqnarray*}
&({\rm I})& {\rm Vol}(FBGH)=0.00628091, \\
&({\rm II})&        {\rm Vol}(FEDG)=0.00486715, \ \ {\rm Vol}(FDGH)=0.00481365, \\
&({\rm III})&       {\rm Vol}(AFED)=0.0189182, \\
&({\rm IV})&       {\rm Vol}(FBGC)=0.0146084, \\
&({\rm V})&        {\rm Vol}(CGFE)=0.00650684, \\
&({\rm VI})&       {\rm Vol}(AFCE)=0.0262516.
\end{eqnarray*}
As an illustration, the calculation
of the volume of $FBGH$ is given in the Appendix.
Note that these values add up to
the volume of the spherical tetrahedron
$T=ABCD
= \{ (x_1, \ldots, x_5) \in \bbS^3
\mid x_1 \le \cdots \le x_5, \ x_{24} \le x_{15} \}:$
\[
{\rm Vol}(T)
=\frac{{\rm Vol}(\bbS^3)}{5!\cdot 2}
=\frac{2\pi^2}{5!\cdot 2}
=0.0822467.
\]

Let
$S=
\{ (x_1, \ldots, x_5) \in \bbS^3
\mid x_1 \le \cdots \le x_5 \}:
{\rm Vol}(S)
=2 {\rm Vol}(T).$
We use the values above to arrive at
\begin{eqnarray*}
{\rm Prob}(({\rm I}) \mid S)
&:=&{\rm Prob}
(x_{14}<x_{23}<x_{24}<x_{15}<x_{25}<x_{34} \mid \\
&& \ \ \ \ \ \ \ \ \ \ \ \ \ \ \ \ \ \ \ \ \ \ \ \
x_1 < \cdots < x_5) \\
&\ =& \frac{{\rm Vol}(FBGH)}{{\rm Vol}(S)} \\
&\ =& \frac{0.00628091}{2\times 0.0822467}=0.0381834.
\end{eqnarray*}

By replacing $x_{i} $ by $-x_{i} $,
we also consider the following cases:

\begin{eqnarray*}
&({\rm I'})&   x_{23}<x_{14}<x_{15}<x_{24}<x_{34}<x_{25}, \\
&({\rm II'})&  x_{14}<x_{23}<x_{15}<x_{24}<x_{34}<x_{25},  \\
&({\rm III'})& x_{14}<x_{15}<x_{23}<x_{24}<x_{34}<x_{25}, \\
&({\rm IV'})& x_{23}<x_{14}<x_{15}<x_{24}<x_{25}<x_{34}, \\
&({\rm V'})&  x_{14}<x_{23}<x_{15}<x_{24}<x_{25}<x_{34}, \\
&({\rm VI'})& x_{14}<x_{15}<x_{23}<x_{24}<x_{25}<x_{34}.
\end{eqnarray*}
Using the symmetry we get
\begin{eqnarray*}
{\rm Prob}(({\rm I}) \mid S) ={\rm Prob}(({\rm I'}) \mid S)
&=&0.0381834,\\
{\rm Prob}(({\rm II}) \mid S)={\rm Prob}(({\rm II'}) \mid S)&=& 0.0588522, \\
{\rm Prob}(({\rm III}) \mid S)={\rm Prob}(({\rm III'}) \mid S)&=& 0.1150086, \\
{\rm Prob}(({\rm IV}) \mid S)={\rm Prob}(({\rm IV'}) \mid S)&=& 0.0888085, \\
{\rm Prob}(({\rm V}) \mid S)={\rm Prob}(({\rm V'}) \mid S)&=& 0.0395569, \\
{\rm Prob}(({\rm VI}) \mid S)={\rm Prob}(({\rm VI'}) \mid S)&=&
0.1595905.
\end{eqnarray*}
We have confirmed that these values coincide with the result of
our simulation study with $\bfx \sim N_5({\bf 0}, I_5),$ 
where $I_5$ denotes the $5 \times 5$-identity matrix.

\section{Concluding remarks}

In this paper, we have solved the problem of counting the number
of ranking patterns in the unidimensional unfolding model
although, due to computational complexity, at present we cannot
determine the explicit number of ranking patterns for $m\ge 9$.
Improving the bound in Lemma \ref{2.4} might reduce the
computational time.  From some computer experiments, it seems that
$L({\cal A}_m)$ and $L({\cal
  A}_{m,q})$ are isomorphic for much smaller $q$ than the value
guaranteed by Lemma \ref{2.4}.

The  problem of counting the number of
ranking patterns can be considered
for the multidimensional unfolding
model. Unlike the unidimensional case, the problem does not reduce to
counting chambers of a hyperplane arrangement and the problem
seems to be quite difficult at this stage.

\section{Appendix}

In this Appendix we illustrate the derivation of the volumes of
the spherical tetrahedra in Section 6 by actually calculating the
volume of $FBGH.$ Take the following orthonormal basis of $H_0'=\{
(x_1,\ldots, x_5)\in \bbR^5 \mid x_1+\cdots + x_5=0 \}:$
\begin{equation*}
\label{basis}
\begin{cases}
\be_1=\frac{1}{\sqrt{2}}(1,-1,0,0,0), \\
\be_2=\frac{1}{\sqrt{6}}(1,1,-2,0,0), \\
\be_3=\frac{1}{\sqrt{12}}(1,1,1,-3,0), \\
\be_4=\frac{1}{\sqrt{20}}(1,1,1,1,-4).
\end{cases}
\end{equation*}
Let $z_1,\ldots, z_4$ be the coordinates
of $\bfx=(x_1,\ldots,x_5) \in H_0'$ in
terms of this basis:
\[
\bfx = z_1\be_1+\cdots+z_4\be_4 =
\begin{pmatrix}
\frac{z_1}{\sqrt{2}}+ \frac{z_2}{\sqrt{6}}
+ \frac{z_3}{\sqrt{12}}+ \frac{z_4}{\sqrt{20}} \\
-\frac{z_1}{\sqrt{2}}+ \frac{z_2}{\sqrt{6}}
+ \frac{z_3}{\sqrt{12}}+ \frac{z_4}{\sqrt{20}} \\
- \frac{2z_2}{\sqrt{6}}+ \frac{z_3}{\sqrt{12}}

+ \frac{z_4}{\sqrt{20}} \\
-\frac{3z_3}{\sqrt{12}}+ \frac{z_4}{\sqrt{20}} \\
-\frac{4z_4}{\sqrt{20}}
\end{pmatrix}.
\]
The conditions $x_{14} \le x_{23}, \ x_{25} \le x_{34}, \ x_3 \le
x_4, \ x_{24} \le x_{15}$ are the binding inequalities of $FBGH$.
They can be written as
\[
\begin{cases}
-\sqrt{3}z_1-z_2+\sqrt{2}z_3 \ge 0, \\
\sqrt{2}z_1-\sqrt{6}z_2-\sqrt{3}z_3+\sqrt{5}z_4 \ge 0, \\
z_2-\sqrt{2}z_3\ge 0, \\
2\sqrt{2}z_1+\sqrt{3}z_3-\sqrt{5}z_4 \ge 0
\end{cases}
\]
in terms of $z_1, \ldots, z_4.$

Now consider the family of
spherical tetrahedra
$T(a), \ 0 \le a \le 1,$ in
$\bbS^3=\{ (z_1, \ldots, z_4) \in \bbR^4 \mid
z_1^2+\cdots+z_4^2=1 \}$ determined by
\begin{equation}
\label{HS}
\begin{cases}
HS_1: \ - a \sqrt{3} z_1-z_2+\sqrt{2}z_3 \ge 0, \\
HS_2: \ \sqrt{2}z_1-\sqrt{6}z_2-\sqrt{3}z_3+\sqrt{5}z_4 \ge 0, \\
HS_3: \ z_2-\sqrt{2}z_3\ge 0, \\
HS_4: \ 2\sqrt{2}z_1+\sqrt{3}z_3-\sqrt{5}z_4 \ge 0.
\end{cases}
\end{equation}
We want to find ${\rm Vol}(T(1)).$

For each $a, \ 0\le a\le 1,$ let $e_{ij}=e_{ij}(a)$ be the edge of
$T(a)$ determined by $HS_i$ and $HS_j, \ 1\le i< j\le 4,$ and
$v_{ijk}=v_{ijk}(a)$ the vertex of $T(a)$ determined by $HS_i,
HS_j$ and $HS_k, \ 1\le i<j<k\le 4.$ Furthermore, denote the
length of $e_{ij}$ and the dihedral angle along $e_{ij}$ by
$\theta_{ij}=\theta_{ij}(a)$ and $\lambda_{ij}=\lambda_{ij}(a), \
1\le i<j\le 4,$ respectively.

The volume ${\rm Vol}(T(a))$ of $T(a)$ can be regarded as a
function of its six dihedral angles $\lambda_{ij},\ 1\le i <j \le
4$. The partial derivatives $\partial{\rm Vol}(T(a))/
\partial \lambda_{ij}, \ 1\le i<j\le 4,$
are given by the following lemma:

\begin{lemma}[Schl\"{a}fli~\cite{sch}]
The partial derivatives of ${\rm Vol}(T(a))$
with respect to $\lambda_{ij}, \ 1 \le i<j\le 4,$
are given as
\[
\frac{\partial{\rm Vol}(T(a))}{\partial
\lambda_{ij}}=\frac{\theta_{ij}}{2}, \
\ \ \ 1\le i<j\le 4.
\]
\end{lemma}

Making use of the above lemma, we will
calculate ${\rm Vol}(T(1))$ in the
following way:
\begin{eqnarray}
\label{int}
{\rm Vol}(T(1))
&=&{\rm Vol}(T(0))
+\int_0^1 \frac{d{\rm Vol}(T(a))}
{da}da \notag \\
&=& {\rm Vol}(T(0))
+\int_0^1 \sum_{1\le i<j\le 4}
\frac{\partial{\rm Vol}(T(a))}
{\partial \lambda_{ij}}
\frac{d\lambda_{ij}}{da}da \\
&=&{\rm Vol}(T(0))
+ \frac{1}{2}\sum_{1\le i<j\le 4}
\int_0^1 \theta_{ij}(a)
\frac{d\lambda_{ij}}{da}da. \notag
\end{eqnarray}

First we calculate $\lambda_{ij}=\lambda_{ij}(a), \ 1\le i<j\le
4.$ By (\ref{HS}) we have
\begin{eqnarray*}
\lambda_{12}
&=&\arccos
\left( \frac{-(-\sqrt{6}a))}{\sqrt{3a^2+3}\sqrt{16}}\right)
=\arccos \left(\frac{\sqrt{2}a}{4
\sqrt{a^2+1}}\right), \\
\lambda_{13}
&=& \arccos \left( \frac{-(-3)}{\sqrt{3a^2+3}\sqrt{3}}\right)
=\arccos \left(\frac{1}{\sqrt{a^2+1}}\right), \\
\lambda_{14}
&=& \arccos
\left( \frac{-(-2\sqrt{6}a+\sqrt{6})}{
\sqrt{3a^2+3}\sqrt{16}}\right)
=\arccos \left(\frac{\sqrt{2}(2a-1)}{4
\sqrt{a^2+1}}\right)
\end{eqnarray*}
and that the other dihedral angles
are constants.

Next we compute the lengths $\theta_{ij}=
\theta_{ij}(a)$ of the edges
$e_{ij}=e_{ij}(a), \ 1\le i<j\le 4.$
These lengths are obtained by
$\theta_{ij}=
\arccos(v_{ijk}\cdot v_{ijl}), \
1\le i<j\le 4, \ k\ne l, \
k,l\notin \{i,j\},$ where $v_{ijk}$ is regarded as
$v_{i'j'k'}$ with $\{ i',j',k' \}=
\{ i,j,k\}, \ 1\le i'<j'<k'\le 4.$
So we begin by finding
the vertices $v_{ijk}=v_{ijk}(a), \ 1\le i<j<k \le 4$ for $0<a<1.$

Vertex $v_{123}$ is obtained by solving
\[
\begin{cases}
 -\sqrt{3}az_1-z_2+\sqrt{2}z_3 = 0, \\
 \sqrt{2}z_1-\sqrt{6}z_2-\sqrt{3}z_3+\sqrt{5}z_4 = 0, \\
 z_2-\sqrt{2}z_3 = 0, \\
 2\sqrt{2}z_1+\sqrt{3}z_3-\sqrt{5}z_4 \ge 0,
\end{cases}
\]
and the other vertices can be
found in a similar manner.
In this way, we obtain
\begin{eqnarray*}
v_{123}&=&(0,-\sqrt{10},-\sqrt{5},-3\sqrt{3})/\sqrt{42}, \\
v_{124}&=&(-\sqrt{10},-\sqrt{30},
-\sqrt{15}(a+1),-3a-7)/\sqrt{8(3a^2+9a+13)}, \\
v_{134}&=&(0,-\sqrt{10},-\sqrt{5},
-\sqrt{3})/\sqrt{18}.
\end{eqnarray*}
Vertex $v_{234}$ is not needed
below.

Using these values of
$v_{123}, v_{124}$ and $v_{134},$ we
get
\begin{eqnarray*}
\theta_{12}&=&
\arccos(v_{123}\cdot v_{124})
=\arccos \left( \frac{7a+18}{
\sqrt{28(3a^2+9a+13)}}\right), \\
\theta_{13}&=&
\arccos(v_{123}\cdot v_{134})
=\arccos \left( \frac{4}{
\sqrt{21}}\right), \\
\theta_{14}&=&
\arccos(v_{124}\cdot v_{134})
=\arccos \left( \frac{\sqrt{3}(4a+11)}{
6\sqrt{3a^2+9a+13}}\right).
\end{eqnarray*}

Finally, we have ${\rm Vol}(T(0))=0,$ since $HS_1$ with $a=0$ and
$HS_3$ in (\ref{HS}) imply $z_2-\sqrt{2}z_3=0.$ Consequently, by
(\ref{int}) and by numerical integration we can evaluate ${\rm
Vol}(T(1))$ as
\begin{eqnarray*}
&& {\rm Vol}(T(1)) \\
&& \ \ \ \ \ =\frac{1}{2}
\int_0^1 \theta_{12}(a)
\frac{d\lambda_{12}}{da}da
+ \frac{1}{2}
\int_0^1 \theta_{13}(a)
\frac{d\lambda_{13}}{da}da
+ \frac{1}{2}
\int_0^1 \theta_{14}(a)
\frac{d\lambda_{14}}{da}da \\
&& \ \ \ \ \ =\frac{1}{2}
\int_0^1 \arccos \left( \frac{7a+18}{
\sqrt{28(3a^2+9a+13)}}\right)
\cdot
\frac{-\sqrt{2}}{(a^2+1)\sqrt{14a^2+16}}
da \\
&& \ \ \ \ \ \ \ \ \ \ \ \ \ + \frac{1}{2}
\arccos \left( \frac{4}{
\sqrt{21}}\right)
\cdot
\left\{ \lambda_{13}(1)-\lambda_{13}(0) \right\} \\
&& \ \ \ \ \ \ \ \ \ \ \ \ \ + \frac{1}{2}
\int_0^1 \arccos \left( \frac{\sqrt{3}(4a+11)}{
6\sqrt{3a^2+9a+13}}\right)
\cdot
\frac{-a-2}{(a^2+1)\sqrt{4a^2+4a+7}}
da \\
&& \ \ \ \ \ =\frac{-0.0810845}{2}
+\frac{\arccos \left( \frac{4}{\sqrt{21}}\right)
\cdot \left( \frac{\pi}{4}-0 \right)}{2}
+\frac{-0.306702}{2} \\
&& \ \ \ \ \ =0.00628091.
\end{eqnarray*}

\newpage
\thispagestyle{empty}
\hspace{-2cm}\includegraphics[width=16cm]{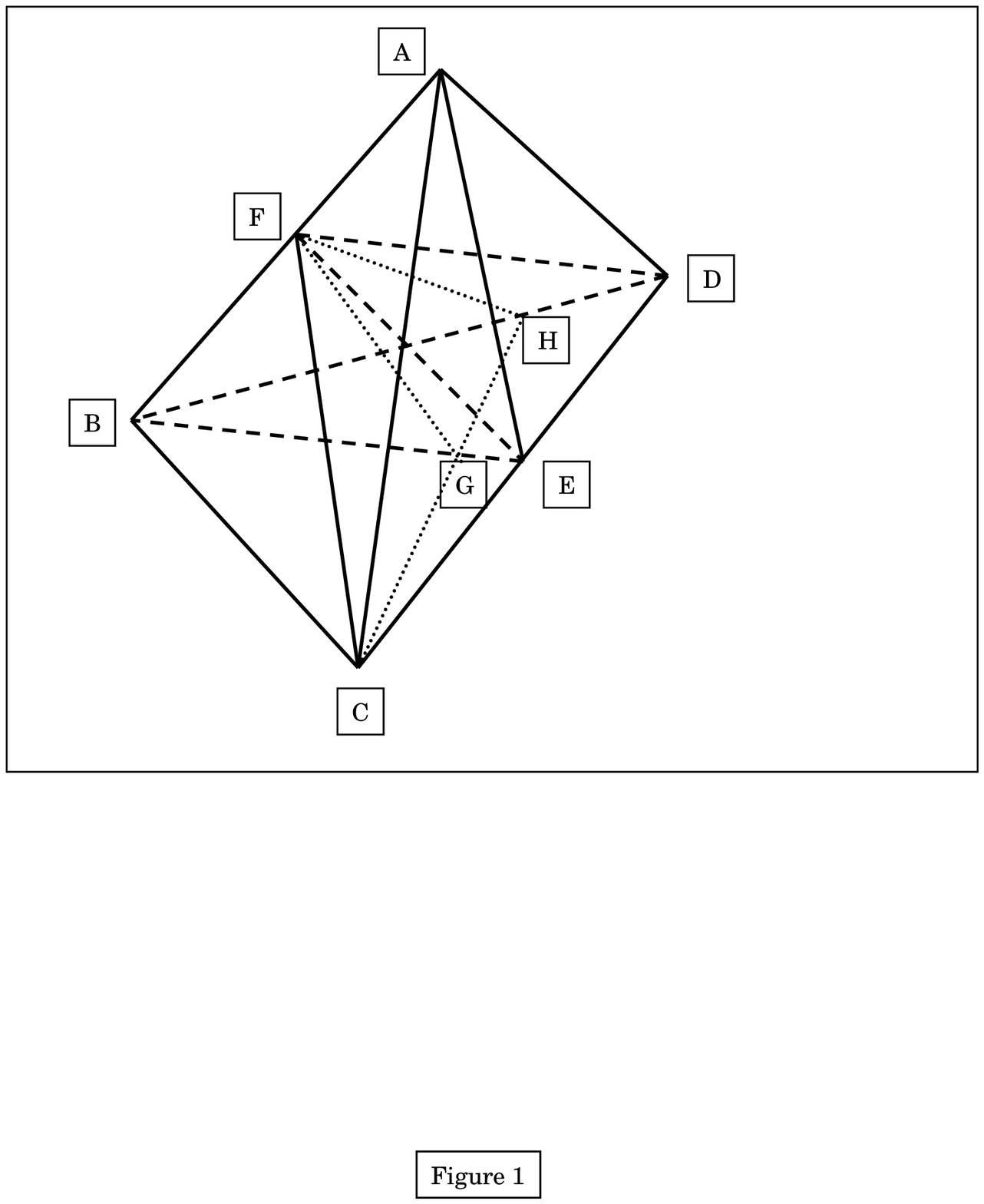}

\newpage
\thispagestyle{empty}
\vspace*{-2cm}
\hspace{-2cm}\includegraphics[width=16cm]{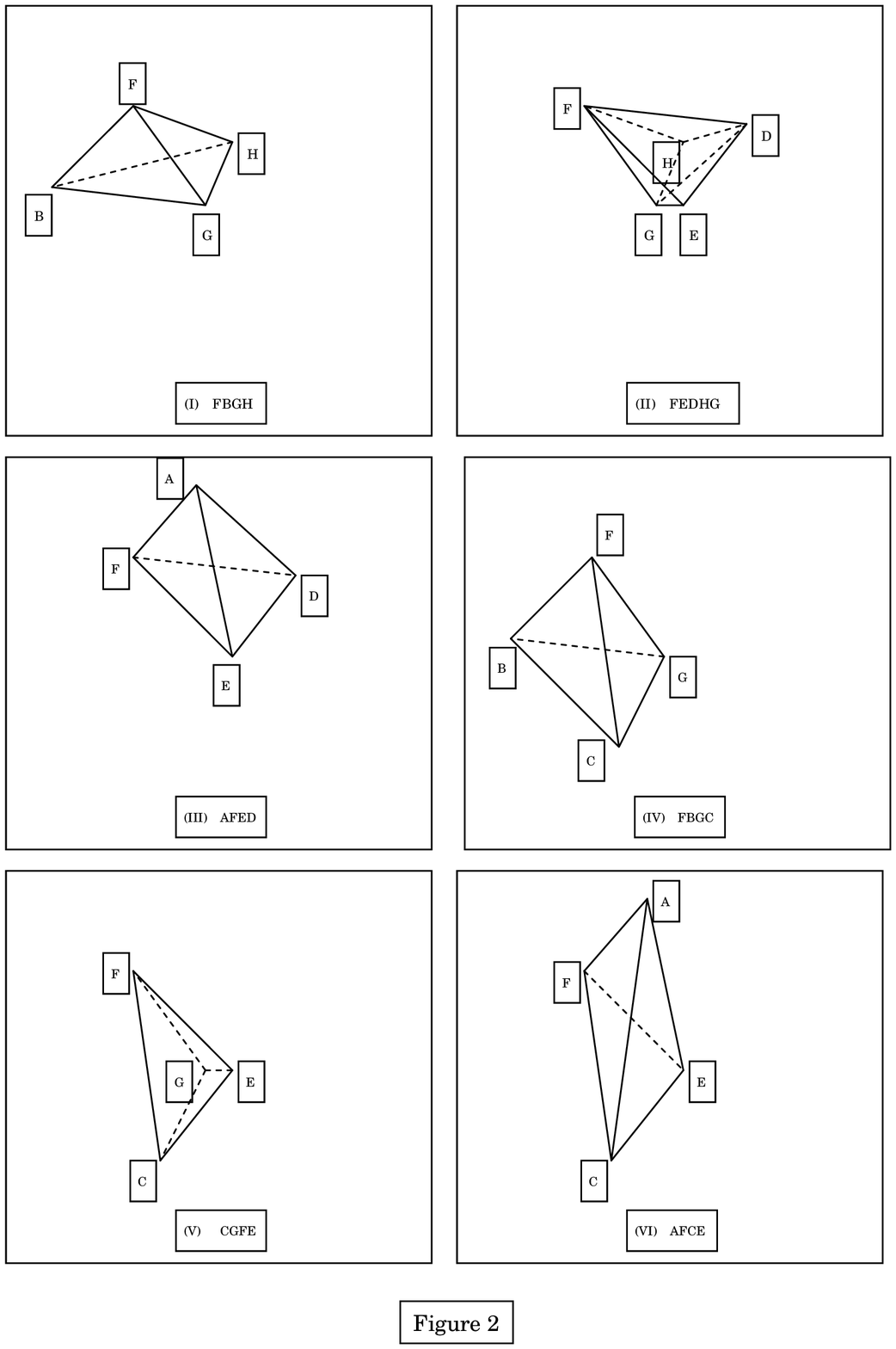}

\end{document}